\documentclass[11pt,reqno,draft]{amsart}

\usepackage[PostScript=dvips,dpi=600,nohug,small]{diagrams}
\usepackage{amssymb}
\usepackage{mathrsfs}
\usepackage{cite}
\usepackage{amsfonts}
\newcommand{\rmv}[1]{}

\def\Span{\mathop{\mathrm{Span}}}
\def\whs{\wh{\chi}_S}
\def\diam#1#2{\mathrm{diam}_{#1}(#2)}
\def\dim{\mathop{\mathrm{dim}}}
\def\QQ{\mathbb Q}
\def\End#1#2{\mathrm{End}_{#1}(#2)}

\def\wh{\widehat}

\def\ZZ{\ensuremath{\mathbb Z}}
\def\RR{\ensuremath{\mathbb R}}

\def\inv{^{-1}}
\def\p{\varphi}

\def\e<{\leq _{E}}

\def\ov#1{\ensuremath{\overline {#1}}}

\def\til#1{\ensuremath{\widetilde {#1}}}

\def\malce{\mathbin{\hbox{$\bigcirc$\rlap{\kern-8.3pt\raise0,50pt\hbox{$\mathtt{m}$}}}}}

\def\CC{\mathbb C}

\def\1sk{^{(1)}}

\def\to{\rightarrow}

\def\Thmname{Theorem}
\def\Propname{Proposition}
\def\Lemmaname{Lemma}
\def\Definitionname{Definition}
\newtheorem{Thm}{\Thmname}[section]
\newtheorem{Prop}[Thm]{\Propname}
\newtheorem{Lemma}[Thm]{\Lemmaname}
{\theoremstyle{definition}
\newtheorem{Def}[Thm]{\Definitionname}}
{\theoremstyle{remark}
\newtheorem{Rmk}[Thm]{Remark}}
\newtheorem{Cor}[Thm]{Corollary}

\newtheorem{Conjecture}{Conjecture}
\newtheorem{Question}{Question}

\numberwithin{equation}{section}

\title{{\v{C}}ern{\'y}'s Conjecture and Group Representation Theory}

\author{Benjamin Steinberg}
\address{School of Mathematics and Statistics\\
Carleton University \\
1125 Colonel By Drive\\
Ottawa, Ontario  K1S 5B6 \\
Canada}
\thanks{The author was supported in part by NSERC}
\email{bsteinbg@math.carleton.ca}
\date{August 10, 2008}

\keywords{{\v{C}}ern{\'y}'s conjecture, synchronizing automata, group representation theory}

\begin{document}
\begin{abstract}
Let us say that a Cayley graph $\Gamma$ of a group $G$ of order $n$ is a {\v{C}}ern{\'y} Cayley graph if every synchronizing automaton containing $\Gamma$ as a subgraph with the same vertex set admits a synchronizing word of length at most $(n-1)^2$.  In this paper we use the representation theory of groups over the rational numbers to obtain a number of new infinite families of {\v{C}}ern{\'y} Cayley graphs.
\end{abstract}
\maketitle

\section{Introduction}
Let $T_X$ be the set of all maps on a set $X$ (which is always taken to be finite in this paper). For the purposes of this article, an \emph{automaton} with state set $X$ is a subset $\Sigma\subseteq T_X$.  Elements of $X$ are commonly referred to as \emph{states}.  Often one writes the automaton as a pair $(X,\Sigma)$ to emphasize the set $X$. Of course, the inclusion $\Sigma\hookrightarrow T_X$ extends to the free monoid $\Sigma^*$ and so an automaton is basically a right action of a finitely generated free monoid on a finite set (where we assume that the generators are sent to different transformations for simplicity).  An important special case is when $G$ is a finite group and $\Delta$ is a generating set for $G$.  The automaton $\Gamma=(G,\Delta)$ is called the \emph{Cayley graph} of $G$ with respect to $\Delta$.  We shall say that an automaton $(G,\Sigma)$ \emph{contains} $\Gamma$ if $\Delta\subseteq \Sigma$.

An important notion in automata theory is that of synchronization.  Let $(X,\Sigma)$ be an automaton.  A word $w\in \Sigma^*$ is called a \emph{synchronizing word} if $|Xw|=1$, that is, $w$ is sent to a constant map under the homomorphism $\Sigma^*\to T_X$.  An automaton that admits a synchronizing word is called a \emph{synchronizing automaton}.  The main open question concerning synchronizing automata is a conjecture from 1964 due to {\v{C}}ern{\'y}~\cite{cerny}, which has received a great deal of attention~\cite{Pincerny,synchgroups,dubuc,cerny,volkovc1,rystsov1,rystsov2,AMSV,trahtman,traht2,volkovc2,Kari,volkovc3}.

\begin{Conjecture}[{\v{C}}ern{\'y}]
A synchronizing automaton with $n$ states admits a synchronizing word of length at most $(n-1)^2$.
\end{Conjecture}

{\v{C}}ern{\'y}, himself, showed that $(n-1)^2$ is the best one can hope for~\cite{cerny}.  The best known upper bound on lengths of synchronizing words is $\frac{n^3-n}{6}$, due to Pin~\cite{twocomb} based on a non-trivial result of Frankl from extremal set theory, see also~\cite{upperbound}.  It should be mentioned that an upper bound of $\frac{n^3-n}{3}$ can be obtained by fairly elementary means, so the hard work lies in improving the bound by a factor of $2$.

There are too many special cases of the {\v{C}}ern{\'y} conjecture that have been proven for us to mention them all here.  The following list of references contain just a few~\cite{Pincerny,synchgroups,dubuc,cerny,volkovc1,rystsov1,rystsov2,AMSV,trahtman,traht2,volkovc2,Kari,volkovc3}.  Let us highlight three results that are most relevant to the paper at hand. We begin with the theorem of Pin~\cite{Pincerny}.

\begin{Thm}[Pin]
Suppose that $\mathscr A=(X,\Sigma)$ is an automaton containing a Cayley graph of a cyclic group of prime order $p$. Then
\begin{enumerate}
\item $\mathscr A$ is synchronizing if and only if some element of $\Sigma$ does not permute $X$;
\item If $\mathscr A$ is synchronizing, then it admits a synchronizing word of length at most $(p-1)^2$.
\end{enumerate}
\end{Thm}

The author (together with his student, Arnold) was motivated by the first part of the above theorem to introduce the notion of a synchronizing group~\cite{synchgroups}: a permutation group $(X,G)$ is said to be a \emph{synchronizing group} if, for each $t\in T_X$ which is not a permutation, the monoid $\langle G,t\rangle$ contains a constant map.  Synchronizing groups have since become a hot topic in the theory of permutation groups~\cite{pneumann} and relate to many classical questions about graphs and finite geometries.  The technique used to study such groups in~\cite{synchgroups} was representation theory over the field of rational numbers, something we explore further in this paper.

Dubuc, in a groundbreaking paper~\cite{dubuc}, extended the second part of Pin's result to Cayley graphs of arbitrary cyclic groups with respect to cyclic generating sets.  This paper was motivated very much by trying to understand Dubuc's ideas from a representation theoretic viewpoint.

\begin{Thm}[Dubuc]
Suppose that $\mathscr A=(X,\Sigma)$ is a synchronizing automaton on $n$ states containing the Cayley graph of a cyclic group with respect to a single generator. Then $\mathscr A$ admits a synchronizing word of length at most $(n-1)^2$.
\end{Thm}

Rystsov~\cite{rystsov1} proved that any synchronizing automaton on $n$ states containing the Cayley graph of a group admits a synchronizing word of length at most $2(n-1)^2$ (this was rediscovered by B\'eal for the special case of cyclic groups~\cite{beal}).  More precisely, if $\Gamma=(G,\Delta)$ is a Cayley graph of a group $G$, the \emph{diameter} $\diam \Delta G$ of $\Gamma$ is the least positive integer $m$ such that each element in $G$ can be represented by an element of $\Delta^*$ of length at most $m$.  Of course $0\leq m\leq |G|-1$.  Rystsov proved the following theorem~\cite{rystsov1}.

\begin{Thm}[Rystsov]\label{rystbound}
Let $\mathscr A=(X,\Sigma)$ be an automaton on $n$ states containing the Cayley graph of a group $G$ with respect to $\Delta$.  Then $\mathscr A$ admits a synchronizing word of length at most $1+(n-1+\diam \Delta G)(n-2)$.
\end{Thm}

Of course, applying the bound of $n-1$ on the diameter yields the upper bound of $2(n-1)^2$.  Notice that Rystsov's bound achieves the {\v{C}}ern{\'y} bound if and only if $\Delta$ contains each non-trivial element of $G$.

If $G$ is a group of order $n$ and $\Delta$ is a set of generators for $G$, then we say that $\Gamma=(G,\Delta)$ is a \emph{{\v{C}}ern{\'y} Cayley} graph if every synchronizing automaton $(G,\Sigma)$ containing $\Gamma$ admits a synchronizing word of length at most \mbox{$(n-1)^2$}.  Let us call $G$ a \emph{{\v{C}}ern{\'y} group} if all its Cayley graphs are {\v{C}}ern{\'y} Cayley graphs.  Of course if the {\v{C}}ern{\'y} conjecture is true, then all groups are {\v{C}}ern{\'y} groups.  Pin's theorem~\cite{Pincerny} establishes that $\ZZ_p$ with $p$ prime is a {\v{C}}ern{\'y} group.  Dubuc~\cite{dubuc} showed that every Cayley graph of $\ZZ_n$ with respect to a cyclic generator is a {\v{C}}ern{\'y} Cayley graph; consequently, $\ZZ_{p^m}$ is a {\v{C}}ern{\'y} group for $p$ prime. To prove that every group is a {\v{C}}ern{\'y} group, one must improve on Rystsov's bound by a factor of $2$.

In this paper, our main result is an improved bound for synchronizing automata containing Cayley graphs based on representation theory over the field of rational numbers.  Our bound does not prove that every Cayley graph is a {\v{C}}ern{\'y} Cayley graph, but it does work for certain Cayley graphs of cyclic groups, dihedral groups, symmetric groups, alternating groups and (projective) special linear groups (in this last example, Galois theory comes into play).  Even when our main result fails to establish a Cayley graph is a {\v{C}}ern{\'y} Cayley graph, our techniques often suffice.  In particular, there are several infinite families of Cayley graphs (coming from affine groups, vector spaces and dihedral groups) that we can prove are {\v{C}}ern{\'y} graphs even though our main result is not up to the task.  As a consequence of our results it follows that if $p$ is a prime, then the dihedral groups $D_p$ and $D_{p^2}$ and the vector spaces $\ZZ_p^m$, for $m\geq 1$, are {\v{C}}ern{\'y} groups.

\section{Representation theory}
As our primary tool in this paper will be representation theory, we try to record here most of the needed background.  There are plenty of excellent books on group representation theory; we shall use~\cite{Dornhoff,curtis} as our primary references.  All groups in what follows should be assumed finite.

\subsection{Basic notions}
Throughout this section $K$ will always be a subfield of the field $\CC$ of complex numbers.
By a representation of a monoid $M$ over $K$, we mean a monoid homomorphism $\p\colon M\to \End K V$ where $\End K V$ is the endomorphism monoid of a finite dimensional $K$-vector space $V$.  It is frequently convenient to denote $\p(m)$ by $\p_m$.   The dimension of $V$ is termed the \emph{degree}  of the representation $\p$, denoted by $\deg(\p)$.  One says that $V$ \emph{carries} or \emph{affords} the representation $\p$.  By the trivial representation of $M$, we mean the homomorphism $\p\colon M\to K=\End K K$ sending all of $M$ to $1$.  If $W\subseteq V$ is a subspace and $A\subseteq M$, we write $AW$ for the subspace spanned by all elements of the form $\p_m(w)$ with $m\in A$ and $w\in W$. A subspace $W\subseteq V$ is said to be \emph{$M$-invariant} if $MW\subseteq W$. Notice that $MW$ is the least $M$-invariant subspace containing $W$.  If $W$ is $M$-invariant, then it affords a \emph{subrepresentation} of $\p$ by restricting each $\p_m$ to $W$.  If the only $M$-invariant subspaces of $V$ are $\{0\}$ and $V$, then $\p$ is said to be \emph{irreducible}.  Evidently every degree one representation is irreducible.

If $\p\colon M\to \End K V$ and $\psi\colon M\to\End K W$ are representations, then their direct sum $\p\oplus \psi\colon M\to \End K{V\oplus W}$ is defined by placing \[(\p\oplus \psi)_m = \p_m\oplus \psi_m.\]
Two representations $\p\colon M\to \End K V$ and $\psi\colon M\to \End K W$ are said to be \emph{isomorphic} (or equivalent) if there is an invertible linear transformation $T\colon V\to W$ such that for each $m\in M$, the diagram
\[\begin{diagram}
V & \rTo^{\p_m} & V\\
\dTo^{T}& & \dTo_T\\
W &\rTo_{\psi_m}& W
\end{diagram}\]
commutes, i.e.\ $\p_m = T\inv \psi_mT$ all $m\in M$.

The \emph{character} $\chi_{\p}$ of a representation $\p$ is the function $\chi_{\p}\colon M\to K$ given by $\chi_{\p}(m) = \mathrm{Tr}(\p_m)$ where $\mathrm{Tr}(A)$ denotes the trace of the linear operator $A$.  Notice that $\chi_{\p}$ only depends on the isomorphism class of $\p$ and $\chi_{\p\oplus \psi} = \chi_{\p}+\chi_{\psi}$.  Also observe that $\chi_{\p}(1)=\deg(\p)$.
A representation $\p$ is said to be \emph{completely reducible} if it is isomorphic to a direct sum of irreducible representations.  The decomposition into irreducibles is unique (up to isomorphism and reordering) and the summands are called the \emph{irreducible constituents} of $\p$.  For a completely reducible representation, every $M$-invariant subspace $W$ has an $M$-invariant complement $W'$ with $V=W\oplus W'$.  Moreover, any irreducible constituent of $V$ is either a constituent of $W$ or of $W'$ (or possibly both if it appears with multiplicity).

It is simple to verify that if \mbox{$\p\colon N\to \End K V$} is an irreducible representation and $\psi\colon M\to N$ is an onto homomorphism, then $\p\psi$ is an irreducible representation of $M$, a fact we shall use without comment.

\subsection{The representation associated to a transformation monoid}
 The primary example of a representation for us is the following.
Let $(X,M)$ be a monoid $M$ acting on the right of a finite set $X$ and let $V=K^X$ be the $K$-vector space of all functions from $X$ to $K$.  Then we can define a representation $\rho\colon M\to \End K V$, called the \emph{standard representation of $(X,M)$},  by right translations: \[\rho_m(f)(x) = f(xm).\] The degree of $\rho$ is, of course $|X|$.  We shall be particularly interested in the case where $M$ is a free monoid, since a pair $(X,\Sigma^*)$ is essentially the same thing as an automaton.
The vector space $V$ comes equipped with the inner product \[\langle f,g\rangle = \sum_{x\in X}f(x)\overline{g(x)}.\]  It is easy to verify that the group of units of $M$ acts by unitary transformations with respect to this inner product.

Let $V_1$ be the subspace of constant functions; let us denote by $\til r$, for $r\in K$, the constant function with value $r$.  Then $\rho_m(\til r) = \til r$, for all $m\in M$, and hence $V_1$ is $M$-invariant.  It is well known and easy to prove that the subspace of vectors fixed by $M$ is precisely the space of constant functions if and only if $M$ acts transitively on $X$.  Set $V_0=V_1^{\perp}$; so $V_0 = \{f\in V: \sum_{x\in X}f(x)=0\}$ and $\dim V_0=|X|-1$.  The subspace $V_0$ is invariant for the group of units of $M$, but not in general for $M$.  It will be convenient to define the \emph{augmentation map} $\epsilon\colon V\to K$ by \[\epsilon(f) = \langle f,\til 1\rangle = \sum_{x\in X}f(x).\]  Observe  that $V_0= \ker \epsilon$.  Let $S\subseteq X$ be a subset and $\chi_S$ its characteristic function.  Then, for $m\in M$, notice $\rho_m(\chi_S) = \chi_{Sm\inv}$ since $\rho_m(\chi_S)(x) = \chi_S(xm)$, which is $1$ if $xm\in S$ and $0$ otherwise.  Also observe that $\epsilon (\chi_S) = |S|$. It is easily verified that
\begin{equation}\label{defhat}
\wh{\chi}_S = \chi_S-\frac{|S|}{|X|}\cdot \til 1=\chi_S-\til{\left(\frac{|S|}{|X|}\right)}
\end{equation}
is the orthogonal projection of $\chi_S$ onto $V_0$.  Indeed, the vector $\til{1}$ spans $V_1$ and $\dfrac{\langle \chi_S,\til{1}\rangle}{\langle \til 1,\til 1\rangle} = \dfrac{|S|}{|X|}$.  Notice that $\whs=0$ if and only if $S=X$.  The following observation will be applied often in this paper.

\begin{Prop}\label{constantout}
Suppose that $(X,M)$ is transitive and let $\rho$ be the standard representation of $(X,M)$.  Assume that some element of $M$ acts as a constant map on $X$.  Let $S$ be a proper subset of $X$ and set $W=\Span\{\whs\}$. Then $MW\nsubseteq V_0$.
\end{Prop}
\begin{proof}
By transitivity of $M$, there is a constant map $f\in M$ with image contained in $S$.  We then compute \[\epsilon(\rho_f(\whs)) = \epsilon(\chi_{Sf\inv}) -|S| = |Sf\inv|-|S|=|X|-|S|>0\] and so $\rho_f(\whs)\notin \ker\epsilon=V_0$.
\end{proof}

\begin{Rmk}
The following remark is for experts in representation theory.
If $M$ acts faithfully and transitively on $X$ and contains a constant map, then one can verify that the standard representation of $(X,M)$ is an injective indecomposable representation with simple socle $V_1$.
\end{Rmk}

A fact that we shall use frequently is that if $G$ is a finite group acting transitively on $X$, then $\frac{1}{|G|}\sum_{g\in G}\rho_g$ is the orthogonal projection of $V$ onto $V_1$ and, in particular, it annihilates $V_0$.

\begin{Prop}\label{orthogonalproj}
Let $G$ be a finite group acting transitively on the right of a finite set $X$.  Then \[P = \frac{1}{|G|}\sum_{g\in G}\rho_g\] is the orthogonal projection onto $V_1$.
\end{Prop}
\begin{proof}
Since $V_0, V_1$ are both $G$-invariant, they are both invariant under $P$.  So if we can show $V_1=\mathop{\mathrm {Im}}P$ and $P$ fixes $V_1$, then the proposition will follow from the orthogonal decomposition $V=V_0\oplus V_1$.  Let us prove the latter statement first.  Since each element of $G$ fixes $V_1$, if $\til r\in V_1$, then \[P\til r = \frac{1}{|G|}\sum_{g\in G}\rho_g(\til r) = \frac{1}{|G|}\sum_{g\in G}\til r = \til r\] and hence $P$ fixes $V_1$.  Next let $f\in V$ and let $x,y\in X$.  By transitivity $y=xh$ some $h\in G$.  Then we have \[Pf(y) = \frac{1}{|G|}\sum_{g\in G}f(yg) = \frac{1}{|G|}\sum_{g\in G}f(xhg)=\frac{1}{|G|}\sum_{t\in G}f(xt) = Pf(x)\] where the last equality follows by making the change of variables $t=hg$.  It follows that $Pf$ is a constant map, completing the proof.
\end{proof}

Since $P$ obviously fixes any vector fixed by all of $G$, the above proposition shows that $V_1$ is the space of fixed vectors of $G$, as was mentioned earlier.

\subsection{Group representation theory}
We highlight here some key points about group representations. Let $G$ be a finite group.  Maschke's theorem says that every representation of $G$ over $K$ is completely reducible~\cite{Dornhoff,curtis}.  It is a standard fact that group representations are determined up to isomorphism by their characters~\cite[Chpt.~7]{Dornhoff} and hence often one does not distinguish between an irreducible representation and its associated character.  If we let $G$ act on the right of itself by right multiplication, then the standard representation of $(G,G)$ is called the \emph{regular representation} of $G$.  It is well known that each irreducible representation (up to isomorphism) of $G$ is a constituent in the regular representation of $G$.  In particular, if we look at the decomposition of the regular representation into $V_0\oplus V_1$, then we see that each non-trivial irreducible representation of $G$ is a constituent of $V_0$ and each constituent of $V_0$ is non-trivial.  Representations $\rho$ and $\psi$ are said to be \emph{orthogonal} if they have no common irreducible constituents.  Then, we have the following consequence of Proposition~\ref{orthogonalproj}.

\begin{Prop}\label{annihirred}
Let $G$ be a group and $\p\colon G\to \End K V$ be a representation of $G$ orthogonal to the trivial representation.  Then \[0=\frac{1}{|G|}\sum_{g\in G}\p_g.\]
\end{Prop}
\begin{proof}
Each irreducible constituent of the representation $\p$ is an irreducible constituent of $V_0$ in the regular representation and hence is annihilated by $\frac{1}{|G|}\sum_{g\in G}\p_g$ thanks to Proposition~\ref{orthogonalproj}.
\end{proof}

If $K=\CC$, then the number of isomorphism classes of irreducible representations of $G$ is precisely the number of conjugacy classes of $G$.  Moreover, if $\p^{(1)},\ldots,\p^{(s)}$ form a complete set of representatives of the equivalence classes of irreducible representations of $G$ over $\CC$ and $d_i$ is the degree of $\p^{(i)}$, then $\p^{(i)}$ appears exactly $d_i$ times as a summand in the decomposition of the regular representation of $G$ into irreducibles.  In particular, $|G|=d_1^2+\cdots+d_s^2$, see~\cite{Dornhoff,curtis}.

Every representation of $G$ over $\QQ$ is isomorphic to a matrix representation $\p\colon G\to M_n(\QQ)$ where $M_n(\QQ)$ is the monoid of $n\times n$-matrices over $\QQ$ (simply choose a basis for the representation space).  Hence each representation over $\QQ$ can be viewed as a representation over $\CC$. (Formally speaking, one replaces $V$ by the tensor product $\CC\otimes_{\QQ} V$.) One says that $\p$ is \emph{absolutely irreducible} if it is irreducible as a representation over $\CC$.  Absolutely irreducible representations must be irreducible, but not conversely.  For example, let $\omega_n$ be a primitive $n^{th}$-root of unity.  Then one can define an irreducible representation $\p\colon \ZZ_n\to \End {\QQ} {\QQ(\omega_n)}$ by having the generator act via left multiplication by $\omega_n$.  It is easy to see that a $\ZZ_n$-invariant subspace is the same thing as a left ideal in $\QQ(\omega_n)$, but $\QQ(\omega_n)$ is a field and so has no non-zero proper ideals.  However, every irreducible representation of $\ZZ_n$ over $\CC$ has degree $1$ (since it has $n$ conjugacy classes and the sums of the degrees squared add up to $n$). So $\p$ is not absolutely irreducible.

It is a classical fact that if $\chi$ is the character of a complex representation of a group $G$ of order $n$, then $\chi(g)$ is a sum of $n^{th}$-roots of unity and hence is an algebraic number (in fact an algebraic integer), for each $g\in G$~\cite{Dornhoff,curtis}.  Thus one can form a number field $\mathbb Q(\chi)$ (i.e.\ a finite extension of $\QQ$), called the \emph{character field of $\chi$}, by adjoining the values of $\chi$.  In fact, $\mathbb Q(\chi)$ is a subfield of the cyclotomic field $\mathbb Q(\omega_n)$ and therefore is a Galois (in fact abelian) extension of $\QQ$.  Hence if \mbox{$H=\mathrm{Gal}(\QQ(\chi):\QQ)$} is the Galois group of $\QQ(\chi)$ over $\QQ$, then $|H|=[\QQ(\chi):\QQ]$.  Notice that $H$ acts on the right of the set of functions $\theta\colon G\to \QQ(\chi)$ by putting $\theta^{h}(g) = h\inv(\theta(g))$ for $h\in H$.  The main result of~\cite[Chpt.~24]{Dornhoff} establishes the following theorem, encapsulating the relationship between irreducible representations of $G$ over $\QQ$ and $\CC$.

\begin{Thm}\label{schurindex}
Let $G$ be a finite group.
\begin{enumerate}
\item  Let $\theta$ be the character of an irreducible representation of $G$ over $\mathbb Q$.  Then there is a complex irreducible character $\chi$ of $G$ and an integer $s(\chi)$, called the Schur index of $\chi$, so that \[\theta = s(\chi)\cdot \sum_{h\in \mathrm{Gal}(\QQ(\chi):\QQ)}\chi^{h}.\]
\item If $\chi$ is the character of complex irreducible representation of $G$, then there is a unique integer $s(\chi)$ so that \[\theta = s(\chi)\cdot \sum_{h\in \mathrm{Gal}(\QQ(\chi):\QQ)}\chi^{h}\] is the character of an irreducible representation of $G$ over $\QQ$. In particular, one has
    \begin{equation}\label{schurestimate}
    \deg(\theta) = s(\chi)[\QQ(\chi):\QQ]\deg(\chi)\geq [\QQ(\chi):\QQ]\deg(\chi).
    \end{equation}
\end{enumerate}
\end{Thm}

Hence the representation theory of $G$ over $\QQ$ can be understood in principle via the complex representation theory and Galois theory.   However, it should be mentioned that computing the Schur index is a non-trivial task and so we content ourselves in this paper with the bound in~\eqref{schurestimate}.

\subsection{Representations of free monoids}
Several combinatorial lemmas concerning representations of free monoids have been exploited in the literature in connection with {\v{C}}ern{\'y}'s conjecture~\cite{Kari,dubuc,beal,rystsov1}, as well as with the theory of rational power series~\cite{BerstelReutenauer}.  Here we collect some variants.  Let us denote by $\Sigma^{\leq d}$ the set of all words in $\Sigma^*$ of length at most $d$.  The length of a word $w$ is denoted $|w|$, as usual.

\begin{Lemma}\label{chainin}
Let $\p\colon \Sigma^*\to \End K V$ be a representation and suppose that $W\subseteq V$ is a subspace.  Then $\Sigma^*W = \Sigma^{\leq d}W$ where $d=\dim \Sigma^*W-\dim W$.
\end{Lemma}
\begin{proof}
Let $W_i=\Sigma^{\leq i}W$.  Then \[W=W_0\subseteq W_1\subseteq\cdots \subseteq \Sigma^*W\] and if $W_i = W_{i+1}$, then $W_i = \Sigma^*W$.  In particular, if we have \[W_0\subsetneq W_1\subsetneq \cdots \subsetneq W_d=\Sigma^*W,\] then $d+\dim W_0\leq \dim \Sigma^*W$ and so $d\leq \dim \Sigma^*W-\dim W$, as required.
\end{proof}

Our next two results are important for when the alphabet is partitioned into two subsets.

\begin{Lemma}\label{chainpartition}
Let $\Sigma = \Delta\cup \Lambda$ and suppose $\p\colon \Sigma^*\to \End K V$ is a representation.  Let $W\subseteq V$ be a $\Delta^*$-invariant subspace that is not $\Sigma^*$ invariant.  Let $U = \Delta^*\Lambda^{\leq 1}W$.  Then $U = \Delta^{\leq d}\Lambda^{\leq 1}W$ where $d=\dim U-\dim W-1$.
\end{Lemma}
\begin{proof}
By assumption, $W' = \Lambda^{\leq 1}W\supsetneq W$.  Hence $\dim W'\geq \dim W+1$.  Applying the previous lemma to $W'$, we may take \[d=\dim U-\dim W'\leq \dim U-(\dim W+1) = \dim U-\dim W-1,\] establishing the lemma.
\end{proof}

\begin{Prop}\label{erase}
Suppose that $\Sigma = \Delta\cup \Lambda$ and let $\delta\colon \Sigma^*\to \Lambda^*$ be the map erasing letters from $\Delta$.  Let $\p\colon \Sigma^*\to \End K V$ be a representation and $W\subseteq V$ a subspace. Define $W_r = \Span\{wW: |\delta(w)|\leq r\}$ and set
\begin{align*}
V_r &= \Span\{wW: |w|\leq \dim W_r-\dim W, |\delta(w)|\leq r\}&r\geq 0 \\
U_r&=\Delta^{\leq d_r}\Lambda^{\leq 1}V_{r-1} &r\geq 1
\end{align*}
where $d_r = \dim W_{r}-\dim W_{r-1}-1$. Suppose $W_s\neq \Sigma^*W$.  Then, $V_0=W_0$ and, for $1\leq r\leq s+1$, we have $U_r=V_r=W_r$.
\end{Prop}
\begin{proof}
As $W_0=\Delta^*W$, Lemma~\ref{chainin} provides the equality $V_0=W_0$.  It follows directly from the definitions that in general $U_r\subseteq V_r\subseteq W_r$.   Suppose by induction that $V_r=W_r$ for $0\leq r\leq s$; we show $U_{r+1}=V_{r+1}=W_{r+1}$.  Indeed, by induction we have  \[W_{r+1} = \Delta^*\Lambda^{\leq 1}W_r=\Delta^{\leq d_{r+1}}\Lambda^{\leq 1}W_r=\Delta^{\leq d_{r+1}}\Lambda^{\leq 1}V_r=U_{r+1}\] where the second equality  is a consequence of  Lemma~\ref{chainpartition} and the penultimate one follows from the induction hypothesis.  This completes the induction.
\end{proof}

Our final lemma concerns the situation where $W\subseteq U$, but $\Sigma^*W\nsubseteq U$.  The question is how long a word does it take to get you out of $U$?  The answer is provided by the next lemma.

\begin{Lemma}\label{getout}
Suppose  $\p\colon \Sigma^*\to \End K V$ is a representation and let $W\subseteq U$ be subspaces of $V$ such that $\Sigma^*W\nsubseteq U$. Let us say $W = \mathop{\mathrm{Span}}X$.  Then there exist $x\in X$ and $w\in \Sigma^*$ with $\p_w(x)\notin U$ and $|w|\leq \dim U-\dim W+1$.
\end{Lemma}
\begin{proof}
Again let $W_i=\Sigma^{\leq i}W$ and consider the chain of subspaces \[W=W_0\subseteq W_1\subseteq\cdots \subseteq \Sigma^*W.\]  As in the proof of Lemma~\ref{chainin}, if $W_i = W_{i+1}$, then $W_i = \Sigma^*W$. Now $W_0\subseteq U$ and $\Sigma^*W\nsubseteq U$, so choosing $d$ least such that $W_d\nsubseteq U$, we have \[W=W_0\subsetneq W_1\subsetneq \cdots \subsetneq W_{d-1}\subseteq U.\]   Consequently, $\dim W+d-1\leq \dim U$, or in other words we have the sought after inequality $d\leq \dim U-\dim W+1$.  Since $W_d$ is spanned by the elements $\p_w(x)$ with $|w|\leq d$ and $x\in X$, and moreover $W_d\nsubseteq U$, it follows that we can find $w\in \Sigma^*$ and $x\in X$ with the desired properties.
\end{proof}

\section{An improved bound for automata containing Cayley graphs}
In this section, we ameliorate Rystsov's bound for synchronizing automata containing Cayley graphs. Our bound is good enough to obtain Pin's result~\cite{Pincerny}, as well as to obtain several new infinite families of {\v{C}}ern{\'y} Cayley graphs.  It does not recover Dubuc's result, although it comes much closer than~\cite{beal,rystsov1}. Again all groups are finite here.

Let $G$ be a group of order $n>1$.  Define $m(G)$ to be the maximum degree of an irreducible representation of $G$ over $\QQ$.  As each irreducible representation of $G$ is a constituent in the regular representation, and all groups admit the trivial representation, one has $1\leq m(G)\leq n-1$.  Since the regular representation is faithful, it follows from Maschke's theorem that the irreducible representations of $G$ separate points.  Since $\QQ^*\cong \ZZ_2$, it follows that $m(G)=1$ if and only if $G\cong \ZZ_2^m$ for some $m$.  We shall see momentarily that if $G$ is a cyclic group of prime order $n$, then $m(G)=n-1$. Before proving our main theorem, we isolate some key ideas of the proof, many of which are inspired by the beautiful paper of Dubuc~\cite{dubuc}; see also our previous paper with Arnold~\cite{synchgroups}.

Let $(X,\Sigma)$ be a synchronizing automaton with $n$ states.  Suppose $\Sigma^*$ acts transitively on $X$ (as is usually the case).  Then for any proper subset $S\subseteq X$, there is a word $w\in \Sigma^*$ so that $|Sw\inv|>|S|$: one can take $w$ to be an appropriate synchronizing word, for instance.
The basic  strategy for obtaining bounds on lengths of synchronizing words (although this strategy is now known not to be optimal in general~\cite{Karicounter}) is to prove that, for each subset $S$ of $X$ with $2\leq|S|<n$, there is a word $t\in \Sigma^*$ of length at most $k$ so that $|St\inv|>|S|$ (we say such a word $t$ \emph{expands} $S$).  Then one obtains a synchronizing word of length at most $1+k(n-2)$.  Indeed, to expand a singleton subset requires a single non-permutation from $\Sigma$ (which must exist if the automaton is synchronizing).  One can then expand repeatedly by words of length at most $k$ until obtaining $X$.  Since one has to expand at most $n-2$ times from a two element set to an $n$ element set, this establishes the bound.  Observing that $(n-1)^2= 1+n(n-2)$, the goal is to try and prove that one can take $k\leq n$.

Our first idea is a lemma that we shall refer to as the ``Standard Argument'' since it is an argument we shall use time and time again throughout the paper.

\begin{Lemma}[Standard Argument]
Suppose $(X,\Sigma)$ is an automaton and let $\rho\colon \Sigma^*\to \End {\QQ} V$ be the standard representation with $V=\QQ^X$. Let $V_1$ be the space of constant maps and $V_0$ the orthogonal complement.  Let $S\subsetneq X$ and recall the definition of $\whs$ from \eqref{defhat}.  Suppose  $\rho_{uvw}(\whs)\notin V_0$ with $u,v,w\in \Sigma^*$.  Then if there exist $r\geq |v|$ and a non-negative linear combination \[P=\sum_{y\in \Sigma^{\leq r}}c_y\rho_y\] with $c_v>0$  and $\rho_uP\rho_w(\whs)\in V_0$, then $|St\inv|>|S|$ for some $t\in \Sigma^*$ with $|t|\leq |u|+|w|+r$.
\end{Lemma}
\begin{proof}
Since $\rho_{uvw}(\whs)\notin V_0=\ker \epsilon$, it follows
\[0\neq \epsilon(\rho_{uvw}(\whs)) = \epsilon(\chi_{S(uvw)\inv})-\epsilon\left(\frac{|S|}{|X|}\cdot \til 1\right) = |S(uvw)\inv| - |S|.\] This leads us to two cases.  If $|S(uvw)\inv| - |S|>0$, then we are done since $|uvw|=|u|+|v|+|w|\leq |u|+|w|+r$.  So suppose instead
\begin{equation}\label{negative}
|S(uvw)\inv| - |S|<0.
\end{equation}
Since $\rho_uP\rho_w(\whs)\in V_0=\ker \epsilon$, it follows
\begin{equation}\label{zerosum}
\begin{split}
0=\epsilon(\rho_uP\rho_w(\whs)) &= \sum_{y\in \Sigma^{\leq r}} c_y\epsilon(\rho_{uyw}(\whs))\\ &= \sum_{y\in \Sigma^{\leq r}}c_y(|S(uyw)\inv| - |S|).
\end{split}
\end{equation}
Taking into account that the $c_y$ are non-negative, $c_v>0$ and \eqref{negative} holds, in order for \eqref{zerosum} to be valid there must exist $y\in \Sigma^{\leq r}$ with $|S(uyw)\inv| - |S|>0$.  Setting $t=uyw$ completes the proof.
\end{proof}

The next lemma, which shall be our main workhorse, is called the ``Gap Bound''.  First let us describe the ``Standard Setup'', which is essentially a collection of notational conventions that will be needed at the start of nearly every proof in the remainder of the paper.

\begin{Def}[Standard Setup]
Let $G$ be a group of order $n>1$ generated by $\Delta$ and suppose $\Delta\subseteq \Sigma\subseteq T_G$ with $(G,\Sigma)$ a synchronizing automaton. Set $\Lambda = \Sigma\setminus \Delta$.  Suppose $S\subseteq G$ is a subset with $2\leq |S|< n$. Let $\rho\colon \Sigma^*\to \End {\QQ} V$ be the standard representation where $V=\QQ^G$.   Put $W = \Span \{\whs\}$ and set $W_r = \Span\{wW: |\delta(w)|\leq r\}$, for $r\geq 0$, and we agree $W_{-1}=0$.  Recall that $\delta\colon \Sigma^*\to \Lambda^*$ is the map erasing $\Delta$.  Define $c_r=\dim W_r-\dim W_{r-1}$.  These numbers are referred to as the \emph{gaps}.  By construction $W_r$ is a $G$-invariant subspace for the regular representation of $G$ so we may write $W_r=W_{r-1}\oplus U_r$ where  $U_r$ is a $G$-invariant subspace. Note that $c_r=\dim U_r$ and $W_r=U_0\oplus U_1\oplus\cdots \oplus U_r$.   Then \[W\subseteq W_0\subseteq W_1\subseteq\cdots \subseteq \Sigma^*W\] and as soon as $W_r=W_{r+1}$ one has $W_r=\Sigma^*W$ (since $W_r=\Delta^*\Lambda^{\leq 1}W_{r-1}$ for $r\geq 1$). Note that $W_0=\Delta^*W\subseteq V_0$, while Proposition~\ref{constantout} yields $\Sigma^*W\nsubseteq V_0$.  Hence there is a maximal integer $s$ so that $W_s\subseteq V_0$.
\end{Def}

The Gap Bound relates the length of a word needed to expand $S$ to the size of the maximal gap.

\begin{Lemma}[Gap Bound]
Let us assume the Standard Setup.  Then there is a word $t\in \Sigma^*$ of length at most \[1+\dim W_s-\max_{0\leq r\leq s}\{c_r\}+\diam \Delta G\] such that $|St\inv|>|S|$.
\end{Lemma}
\begin{proof}
Fix, for each element $g\in G$, a word $u_g\in \Delta^*$ of length at most $\diam \Delta G$ so that $u_g$ maps to $g\in G$ under the projection \mbox{$\pi\colon \Delta^*\to G$}.  Let $\lambda\colon G\to \End {\QQ} V$ be the regular representation of $G$. Then $\rho|_{\Delta^*}=\lambda \pi$, that is, $\rho_{u} = \lambda_{\pi(u)}$ for $u\in \Delta^*$.  In particular, $\rho_{u_g} = \lambda_g$.

Since $W_{s+1}\nsubseteq V_0$, it follows $\Lambda W_s\nsubseteq V_0$ as $V_0$ is invariant under $\Delta^*$. Hence $bW_s\nsubseteq V_0$ some $b\in \Lambda$.  Let $c_k=\max\{c_r:0\leq r\leq s\}$.
First suppose $k=0$.  Proposition~\ref{erase}, but with $W_0$ in the place of $W$, implies that $W_s$ is spanned by elements of the form $\rho_x(f)$ where $|x|\leq \dim W_s-\dim W_0=\dim W_s-c_0$, $|\delta(x)|\leq s$ and $f\in W_0$. As $W_0=\Delta^*W$, it follows $\rho_{bxy}(\whs)\notin V_0$ for some $y\in \Delta^*$ and $x$ as above.   Since $\whs\in V_0$, we have by Proposition~\ref{orthogonalproj}
\[0=\rho_{bx}\frac{1}{|G|}\sum_{g\in G}\lambda_g(\whs).\]
Recalling that $\rho_y = \lambda_{\pi(y)}=\rho_{u_{\pi(y)}}$, the Standard Argument with $u=bx$, $v=u_{\pi(y)}$, $w=1$ and $P=\frac{1}{|G|}\sum_{g\in G}\rho_{u_g}$  provides a word $t$ of length at most \[|bx|+\diam \Delta G\leq 1+\dim W_s-c_0+\diam \Delta G\] such that $|St\inv|>|S|$.

Finally suppose  $k>0$.  Proposition~\ref{erase} yields $W_k$ is spanned by vectors of the form $\rho_{yb'z}(\whs)$ where $y\in \Delta^*$, $b'\in \Lambda^{\leq 1}$ and $z\in \Sigma^*$ such that the inequalities $|z|\leq \dim W_{k-1}-1$ and $|\delta(z)|\leq k-1$ hold.  On the other hand, Proposition~\ref{erase}, but with $W_k$ in the place of $W$, yields that $W_s$ is spanned by elements of the form $\rho_x(f)$ where $|x|\leq \dim W_s-\dim W_k$, $|\delta(x)|\leq s-k$ and $f\in W_k$.  Putting this together, we can find $x,y,b',z$ with the above properties so that $\rho_{bxyb'z}(\whs)\notin V_0$.   Since $\rho_{b'z}(\whs)\in W_k\subseteq V_0$, it follows from Proposition~\ref{orthogonalproj} that
\[0=\rho_{bx}\frac{1}{|G|}\sum_{g\in G}\lambda_g\rho_{b'z}(\whs).\]
Invoking the Standard Argument where we take $u=bx$, $v=u_{\pi(y)}$, $w=b'z$ and $P=\frac{1}{|G|}\sum_{g\in G}\rho_{u_g}$ yields the existence of a word $t\in \Sigma^*$ with
\begin{align*}
|t|&\leq |bx|+|b'z|+\diam \Delta G\\ &\leq 1+\dim W_s-\dim W_k+1+\dim W_{k-1}-1+\diam \Delta G\\ & = 1+\dim W_s-c_k+\diam \Delta G\end{align*} such $|St\inv|>|S|$.  This completes the proof.
\end{proof}


Since the largest gap is at least $m(G)$, or $n-1-\dim W_s\geq m(G)$, we obtain our main result, improving upon Rystsov's bound, Theorem~\ref{rystbound}.

\begin{Thm}\label{mainthm}
Let $G$ be a group of order $n>1$ generated by $\Delta$ and suppose $\Delta\subseteq \Sigma\subseteq T_G$ with $(G,\Sigma)$ a synchronizing automaton.  Then $(G,\Sigma)$ admits a synchronizing word of length at most \[1+\left(n-m(G)+\diam \Delta G\right)(n-2).\]  In particular, if $\diam \Delta G\leq m(G)$, then $(G,\Sigma)$ satisfies the {\v{C}}ern{\'y} bound and hence $(G,\Delta)$ is a {\v{C}}ern{\'y} Cayley graph.
\end{Thm}
\begin{proof}
Observing that $(n-1)^2 = 1+n(n-2)$, the last statement follows from the first, which we proceed to prove. Let $S\subseteq G$ be a subset with $2\leq |S|< n$.  It suffices to show that there exists $t\in \Sigma^*$ with $|t|\leq n-m(G)+\diam \Delta G$ and $|St\inv|>|S|$.  So we assume the Standard Setup. Let $\theta$ be an irreducible character of $G$ of degree $m(G)$.  We know that $\theta$ appears as a constituent in the regular representation of $G$.  As $G$ is non-trivial, we may assume that $\theta$ is not the character of the trivial representation.  Since in the direct sum decomposition $V=V_0\oplus V_1$, the representation afforded by $V_1$ is the trivial representation, it follows that $\theta$ is a constituent in the subrepresentation afforded by $V_0$.   Now we may write $V_0=W_s\oplus U$ with $U$ a $G$-invariant subspace.  Then we have, following the notation of the Standard Setup, \[V_0= W_s\oplus U=U_0\oplus U_1\oplus\cdots \oplus U_s\oplus U.\]  There are two cases.  Suppose first $\theta$ is a constituent of $U_k$, some $0\leq k\leq s$. Then $c_k\geq m(G)$  and therefore it follows
\begin{align*}
1+\dim W_s-\max_{0\leq r\leq s}\{c_r\}+\diam \Delta G& \leq n-c_k+\diam \Delta G\\ &\leq n-m(G)+\diam \Delta G.
\end{align*}
On the other hand, since $n-1=\dim V_0=\dim W_s+\dim U$, if $\theta$ is a constituent of $U$ then $\dim W_s\leq n-1-m(G)$, yielding \[1+\dim W_s-\max_{0\leq r\leq s}\{c_r\}+\diam \Delta G\leq n-m(G)+\diam \Delta G.\]  The desired word $t$ is now provided by the Gap Bound.
\end{proof}

\section{Some examples}
In this section, we consider several natural families of Cayley graphs and determine whether or not we achieve the {\v{C}}ern{\'y} bound with our bound, and if not we see by how much we fail.  In what follows, $\omega_m$ always denotes a primitive $m^{th}$-root of unity.

\subsection{Cyclic groups}
Our first family of examples consists of cyclic groups.
Let $G$ be a cyclic group of order $n$.  Then the regular representation of $G$ is isomorphic to the representation $\rho\colon G\to \End {\QQ} {\QQ[x]/(x^n-1)}$ which sends the generator to left multiplication by $x$.  One has the direct sum decomposition $\QQ[x]/(x^n-1) = \bigoplus_{d\mid n} \QQ(\omega_d)$ where the generator acts on $\QQ(\omega_d)$ via left multiplication by $\omega_d$. An invariant subspace of $\QQ(\omega_d)$ is the same thing as a left ideal, and hence each $\QQ(\omega_d)$ carries an irreducible subrepresentation.  We conclude $m(G) = \phi(n)$, where $\phi$ is Euler's totient function.  If we choose a cyclic generator for $G$, then the diameter of the resulting Cayley graph is $n-1$.  Theorem~\ref{mainthm} thus yields an upper bound of \mbox{$1+(2n-1-\phi(n))(n-2)$} on the length of a synchronizing word.  If $n$ is prime, then $\phi(n)=n-1$ and so we achieve the {\v{C}}ern{\'y} bound, yielding a new proof of Pin's theorem.  In general, we do not obtain Dubuc's result, although we are much closer than~\cite{rystsov1,beal}.  For instance, suppose $n=p^m$ with $p$ prime.  Then one can compute that the ratio of the {\v{C}}ern{\'y} bound to our bound is approximately $1-\frac{1}{p}$ and so is nearly $1$ when $p$ is very large.

On the other hand, suppose $n=p_1\cdots p_k$ is the prime factorization of a square-free number $n$.  Let us consider the natural generating set for $G$ corresponding to the direct product decomposition $G\cong \ZZ_{p_1}\times \cdots\times \ZZ_{p_k}$.  The diameter with respect to this generating set is $(p_1-1)+\cdots +(p_k - 1)$.  On the other hand $m(G) = \phi(n) = (p_1-1)\cdots(p_k-1)$.  It is easy to see that as long as $n$ is odd or $k\geq 3$, one has $(p_1-1)\cdots(p_k-1)\geq (p_1-1)+\cdots +(p_k - 1)$ and so this Cayley graph of $G$ is a {\v{C}}ern{\'y} Cayley graph, something which is not a consequence of the results of~\cite{dubuc}.

\subsection{Dihedral groups}\label{ss:dihedralgroups}
Let $G=D_n$ be the dihedral group of order $2n$.  Let $s$ be a reflection and $r$ be a rotation of order $n$.  Then every element of $D_n$ can be written in one of the forms $sr^k$, $r^ks$ with $k\leq  \left\lfloor \frac{n}{2}\right\rfloor$, $sr^ks$ with $k< \left\lfloor \frac{n}{2}\right\rfloor$ or $r^k$ with $k\leq \left\lceil \frac{n+1}{2}\right\rceil$.  Hence the diameter of $D_n$ with respect to this generating set is at most $\left\lceil \frac{n+1}{2}\right\rceil$.  One can show that $m(D_n) = \phi(n)$.  Let us just establish that $m(D_n)\geq \phi(n)$.  Indeed, define a representation $\rho\colon D_n\to \End {\QQ} {\QQ(\omega_n)}$ by having $\rho(s)$ act via complex conjugation and $\rho(r)$ act via multiplication by $\omega_n$.  We already know this representation is irreducible when restricted to $\langle r\rangle$ and so it is irreducible for $D_n$.

Suppose first that $n=p^k$ with $p$ an odd prime.  Then since $1-1/p\geq 2/3$, the formula $\phi(n)=n(1-\frac{1}{p})$ yields \[\phi(n)-\frac{n+1}{2}\geq \frac{2n}{3}-\frac{n+1}{2} = \frac{n-3}{6}\geq 0\]
since $n\geq 3$. Thus the Cayley graph of $D_n$ with respect to $r,s$ is {\v{C}}ern{\'y}.

Next consider the case $n=p^kq^\ell$ with $p<q$ odd primes.  We claim again that the Cayley graph of $D_n$ with respect to $r,s$ is {\v{C}}ern{\'y}. Indeed, since $1-1/p\geq 2/3$ and $1-1/q\geq 4/5$, from $\phi(n)=n(1-\frac{1}{p})(1-\frac{1}{q})$ it follows
\begin{equation*}
\phi(n) - \frac{n+1}{2}\geq \frac{8n}{15}-\frac{n+1}{2} = \frac{n-15}{30}\geq 0
\end{equation*}
where the last equality uses $n\geq 15$.

The reader should verify that for all other $n$, our bound does not achieve the {\v{C}}ern{\'y} bound.  The bound we obtain is $1+(n-\phi(n)+\left\lceil \frac{n+1}{2}\right\rceil)(n-2)$, which for many $n$ is not far from the {\v{C}}ern{\'y} bound.  For example, for $n=2^m$, one has $\phi(n) = n/2$.  Thus our main result implies that any synchronizing automaton containing the Cayley graph of $D_n$ with respect to $r,s$ has a synchronizing word of length at most $1+(n+1)(n-2) = (n-1)^2+n-2$.

We shall establish later that if $p$ is an odd prime, then $D_p$ and $D_{p^2}$ are {\v{C}}ern{\'y} groups.

\subsection{Symmetric and alternating groups}
It is well known that each irreducible representation of the symmetric group $S_n$ over $\QQ$ is absolutely irreducible~\cite{curtis}.  Letting $p_n$ be the number of partitions of $n$, it follows that $S_n$ has $p_n$ irreducible representations over $\QQ$ and the sum of their degrees squared is $n!$.  Thus $p_nm(S_n)^2\geq n!$ and so we obtain $m(S_n)\geq \sqrt{n!/p_n}$.  It is a well-known result of Hardy and Ramanujan that $p_n\sim \frac{\exp\left(\pi\sqrt{2n/3}\right)}{4n\sqrt{3}}$.  On the other hand, Stirling's formula says that $n!\sim \sqrt{2\pi n}\left(\frac{n}{e}\right)^n$.  Comparing these expressions, we see that $m(S_n)$ grows faster than any exponential function of $n$.  On the other hand, the diameter of $S_n$ with respect to any of its usual generating sets grows polynomially with $n$. For instance, if one uses the Coxeter-Moore generators $(12),(23),\ldots, (n-1n)$ the diameter of $S_n$ is well known to be $\binom{n}{2}$, while if one uses the generators $(12),(12\cdots n)$, then the diameter is no bigger than $(n+1)n(n-1)/2$ since each Coxeter-Moore generator can be expressed as a product of length at most $n+1$ in these generators.  Thus the Cayley graph of $S_n$ with respect to either of these generating sets is a {\v{C}}ern{\'y} Cayley graph for $n$ sufficiently big (and sufficiently big is not very big in this case).

To deal with the alternating group $A_n$, we use the following lemma, which is a trivial consequence of Clifford's theorem.

\begin{Lemma}\label{index2}
Let $G$ be a group and suppose $H$ is a subgroup of index $2$.  Then $m(H)\geq m(G)/2$.
\end{Lemma}
\begin{proof}
Let $\p\colon G\to \End {\QQ} V$ be an irreducible representation of degree $m(G)$ and fix $s\notin H$.  If $\p|_H$ is irreducible, we are done.  Otherwise, let $W$ be a proper $H$-invariant subspace of $V$ affording an irreducible subrepresentation.  Since $G=H\cup sH$ and $W$ is $H$-invariant, but not $G$-invariant, it follows that $sW\neq W$.   Moreover, $sW$ is also an $H$-invariant subspace since if $h\in H$ and $w\in W$, then $hsw = s(s\inv hs)w\in sW$ using that $H$ is a normal subgroup and $W$ is $H$-invariant.  Moreover, $sW$ carries an irreducible subrepresentation of $H$ since if $U\leq sW$ is an $H$-invariant subspace, a routine verification yields $s\inv U$ is an $H$-invariant subspace of $W$.  Consequently $W\cap sW=0$.  Clearly the direct sum $W\oplus sW$ is $G$-invariant, being preserved by both $H$ and $s$ and using $G = H\cup sH$.  Thus, because $\p$ is irreducible, we conclude that $V=W\oplus sW$.  Since $W$ and $sW$ are isomorphic as vector spaces, $m(G)=\dim V=2\dim W$, establishing the lemma.
\end{proof}

It is immediate from the lemma that $m(A_n)\geq m(S_n)/2$ and hence grows faster than any exponential function of $n$.  Again most of the standard generating sets for $A_n$ have polynomial diameter growth as a function of $n$, leading to {\v{C}}ern{\'y} Cayley graphs for $n$ large enough.

\subsection{Special and projective special linear groups}
Suppose $p$ is an odd prime and let $G=SL(2,p)$ be the group of all $2\times 2$ matrices of determinant $1$ over $\ZZ_p$.  A standard generating set $\Delta$ for $G$ consists of the matrices
\begin{equation}\label{sl2gen}
x=\begin{bmatrix} 1& 1\\ 0 &1\end{bmatrix}\ \text{and}\ y=\begin{bmatrix} 1 &0\\1&1\end{bmatrix}.
\end{equation}
Our goal is to show that the Cayley graph $\Gamma$ of $G$ with respect to $x$ and $y$ is a {\v{C}}ern{\'y} Cayley graph for almost all odd primes.  This is the first example where we shall use the Galois theoretic description of the irreducible representations over $\QQ$.  Let us begin by estimating the diameter, following~\cite{ramanujan}.

A routine computation using $ad-bc=1$ establishes that if $c\neq 0$, then
\begin{equation}\label{diamsl1}
\begin{bmatrix} a & b \\ c& d\end{bmatrix} = \begin{bmatrix} 1 & \frac{a-1}{c} \\ 0& 1\end{bmatrix}\begin{bmatrix} 1 & 0 \\ c& 1\end{bmatrix}\begin{bmatrix} 1 & \frac{d-1}{c} \\ 0& 1\end{bmatrix}.
\end{equation}
On the other hand if $c=0$, then $d\neq 0$ and
\begin{equation}\label{diamsl2}
\begin{bmatrix} a & b \\ 0& d\end{bmatrix} = \begin{bmatrix} a-b & b \\ -d& d\end{bmatrix}\begin{bmatrix} 1 & 0 \\ 1& 1\end{bmatrix}.
\end{equation}
Putting together \eqref{diamsl1} and \eqref{diamsl2} (and using $d\neq 0$ in \eqref{diamsl2} to apply \eqref{diamsl1} to the first matrix in the product) we conclude the diameter $\diam \Delta G$ is at most $3(p-1)+1=3p-2$.

We shall require a lemma about cyclotomic fields for the proof.

\begin{Lemma}\label{adjoincosine}
Let $\alpha = \cos 2\pi/n$ with $n\geq 3$.  Then $[\QQ(\alpha):\QQ]=\phi(n)/2$.
\end{Lemma}
\begin{proof}
The intersection $F$ of $\QQ(\omega_n)$ with the reals $\mathbb R$ is the fixed-field of the automorphism $\sigma\in \mathrm{Gal}(\QQ(\omega_n):\QQ)$ given by $\sigma(z) = \overline z$ (complex conjugation).  Moreover, $\sigma$ is non-trivial as $n\geq 3$ implies $\omega_n\notin \RR$.  Since $\QQ(\omega_n)$ is a Galois extension of $\QQ$, it follows that $[\QQ(\omega_n):F] = |\langle \sigma\rangle| =2$.  Thus
\[\phi(n)=[\QQ(\omega_n):\QQ]=[\QQ(\omega_n):F][F:\QQ]=2[F:\QQ]\] and so $[F:\QQ] = \phi(n)/2$.  Therefore, it suffices to prove $F=\QQ(\alpha)$.  Clearly $\alpha = \frac{1}{2}(\omega_n+\ov{\omega_n})\in F$, so we are left with establishing the containment $F\subseteq \QQ(\alpha)$.  It is easy to see that $\frac{1}{2}(1+\sigma)$ is the projection from $\QQ(\omega_n)$ to $F$ and so $F$ is spanned by the elements $\frac{1}{2}(\omega_n^m+\ov{\omega_n^m}) = \cos 2\pi m/n$ with \mbox{$0\leq m\leq \phi(n)-1$}. Let $T_m$ be the $m^{th}$-Chebyshev polynomial of the first kind~\cite{ramanujan}.  Then $T_m(\cos \theta) = \cos m\theta$.  It follows that $\cos 2\pi m/n$ is a polynomial in $\cos 2\pi/n=\alpha$ and so $F\subseteq \QQ(\alpha)$, as required.
\end{proof}

To conclude the proof, we use the character table of $SL(2,p)$, which goes back to Frobenius and Schur.  It can be found for instance in~\cite[Chpt.~38]{Dornhoff}.  It turns out that $SL(2,p)$ has irreducible complex characters $\chi_1$ of degree $p+1$ with $\QQ(\chi_1) = \QQ(\cos \frac{2\pi}{p-1})$ and $\chi_2$ of degree $p-1$ with character field $\QQ(\chi_2) = \QQ(\cos \frac{2\pi}{p+1})$.  We deduce from Lemma~\ref{adjoincosine} and the estimate \eqref{schurestimate} from Theorem~\ref{schurindex} that
\begin{equation}\label{sl2bound}
m(SL(2,p))\geq \max\left\{(p+1)\frac{\phi(p-1)}{2},(p-1)\frac{\phi(p+1)}{2}\right\}.
\end{equation}
To compare the diameter to $m(SL(2,p))$, first note that $\phi(n)\geq 8$ for all $n>18$. Consequently when our prime $p$ is at least $19$, then \[m(SL(2,p))\geq (p-1)\dfrac{\phi(p+1)}{2}\geq 4(p-1)\geq 3(p-1)+1\] and hence we have a {\v{C}}ern{\'y} Cayley graph.  For $p=17$, a direct computation using \eqref{sl2bound} shows that the graph $\Gamma$ is a {\v{C}}ern{\'y} Cayley graph.  For $p=3,5,7,11,13$ our estimates do not suffice to prove that the graph $\Gamma$ is a {\v{C}}ern{\'y} Cayley graph.

Let us next consider the case of the projective special linear group $G=PSL(2,p) = SL(2,p)/\{\pm I\}$.  We choose the cosets of the matrices $x$ and $y$ from \eqref{sl2gen} as generators and with respect to this generating set, the Cayley graph $\Gamma$ of $G$ still has diameter at most $3(p-1)+1$.  The complex characters of $PSL(2,p)$ are also computed in~\cite[Chpt.~38]{Dornhoff}.  Here one finds an irreducible  character of degree $p+1$ with character field $\QQ(\cos \frac{2\pi}{(p-1)/2})$ and one of degree $p-1$ with character field $\QQ(\cos \frac{2\pi}{(p+1)/2})$.  Arguing as above yields
\begin{equation}\label{psl2}
m(PSL(2,p))\geq \max\left\{\frac{p+1}{2}\cdot\phi\left(\frac{p-1}{2}\right),\frac{p-1}{2}\cdot\phi\left(\frac{p+1}{2}\right)\right\}.
\end{equation}
Again using that $\phi(n)\geq 8$ whenever $n>18$, we conclude that as long as $p\geq 37$, the graph $\Gamma$ is a {\v{C}}ern{\'y} Cayley graph.   Direct computation with the estimate \eqref{psl2} shows that, for $p=19,23,29,31$, we also obtain a {\v{C}}ern{\'y} Cayley graph.  That is, for $p\geq 19$, the Cayley graph of $PSL(2,p)$ with the above generating set is a {\v{C}}ern{\'y} Cayley graph.  Our estimates fail to handle the cases $p=3,5,7,11,13,17$.

\section{Further examples of {\v{C}}ern{\'y} Cayley graphs and {\v{C}}ern{\'y} groups}
In this section we consider some Cayley graphs for which Theorem~\ref{mainthm} is not strong enough to prove that they are {\v{C}}ern{\'y}, but the ideas underlying the theorem do suffice.  In the process we give the first examples of non-cyclic {\v{C}}ern{\'y} groups.  Our main tool is the following lemma, whose proof is similar to that of the Gap Bound.

\begin{Lemma}\label{killwithcyclic}
Assume the Standard Setup.  Let $A$ be a subgroup of $G$ and suppose that, for some $0\leq k\leq s$, one has the decomposition \mbox{$W_k=W_{k-1}\oplus U_k$} where the subspace $U_k$ affords a representation \mbox{$\psi\colon A\to \End {\QQ} {U_k}$} of $A$ so that: each coset of $H=A/\ker \psi$ has a representative in $\Delta^*$ of length at most $c_k$ and either $\psi$ is a non-trivial irreducible representation of $A$, or $A=G$.  Then there exists a word $t$ of length at most $n$ so that $|St\inv|>|S|$.
\end{Lemma}
\begin{proof}
Set $K=\ker \psi$ and choose, for each coset $a\in A/K$, a word \mbox{$u_a\in \Delta^*$} of length at most $c_k$ so that the element of $G$ represented by $u_a$ maps into the coset $a$; without loss of generality, we may assume $u_K=1$.  Let $\Upsilon = \{u_a: a\in A/K\}$.  We view $\psi$ as a representation of $H=A/K$ in the natural way.
First suppose that  $k=0$. Then $W_0=U_0$ and so $W_0$ affords a representation of $H$.  If $A=G$, clearly $HW=GW=W_0$. If $\psi$ is irreducible, then the subrepresentation of $A$ afforded by $W_0$ is irreducible and so again $HW=AW=W_0$.  Applying Lemma~\ref{getout} we can find $u\in \Sigma^*$ with $|u|\leq \dim V_0-\dim W_0+1 = n-c_0$ and $g\in H$ so that $\rho_{uu_{g}}(\whs)\notin V_0$.   Since $W_0$ is orthogonal to the trivial representation of $H$, Proposition~\ref{annihirred} implies $\sum_{a\in H}\psi(a)W_0 = 0$. Thus \[\rho_u\sum_{u_a\in \Upsilon}\rho_{u_a}(\whs)=0.\]  Applying the Standard Argument with $v=u_g, w=1$ yields a word $t$ with $|St\inv|>|S|$ and $|t|\leq |u|+c_0\leq n$.

Next suppose $1\leq k\leq s$.  Then $W_k=W_{k-1}\oplus U_k$ as in the hypothesis. If $\psi$ is irreducible, then since $\Lambda^{\leq 1}W_{k-1}\supsetneq W_{k-1}$ and $W_k/W_{k-1}$ affords an irreducible representation of $A$ isomorphic to $\psi$, factoring by $W_{k-1}$ yields \[W_k/W_{k-1} =  H\Lambda^{\leq 1}W_{k-1}/W_{k-1}.\] It follows $W_k = \Upsilon\Lambda^{\leq 1}W_{k-1}$ (using $1\in \Upsilon$).  On the other hand if $A=G$, then since $W_k = G\Lambda^{\leq 1}W_{k-1}$, it follows that \[W_k/W_{k-1} = G\Lambda^{\leq 1}W_{k-1}/W_{k-1} = H\Lambda^{\leq 1}W_{k-1}/W_{k-1}\] as $W_k/W_{k-1}$ affords a representation isomorphic to $\psi$ and $H = G/\ker \psi$.  So again we have $W_k = \Upsilon\Lambda^{\leq 1}W_{k-1}$.
Now by choice of $s$, we have \mbox{$\Lambda W_s\nsubseteq V_0$}.  Applying Proposition~\ref{erase} with $W_k$ in place of $W$ and $W_s$ in place of $W_r$ it follows that $W_s$ is spanned by vectors of the form $\rho_u(f)$ so that \mbox{$|\delta(u)|\leq s-k$}, $|u|\leq \dim W_s-\dim W_k$ and $f\in W_k$. Hence we can find $b\in \Lambda$ and $u,f$ as above with $\rho_{bu}(f)\notin V_0$.  Now from $W_k = \Upsilon\Lambda^{\leq 1}W_{k-1}$, it follows that we may find such an $f$ of the form $\rho_{u_gb'w}(\whs)$ with $g\in H$, $b\in \Lambda^{\leq 1}$, $|\delta(w)|\leq k-1$ and $|w|\leq \dim W_{k-1}-1$ (again using Proposition~\ref{erase}).  The operator $P=\sum_{u_a\in \Upsilon}\rho_{u_a}$ annihilates $U_k$ by Proposition~\ref{annihirred} (since $U_k$ affords a representation of $H$ orthogonal to the trivial representation) and therefore $PW_k\subseteq W_{k-1}$.  Since $\rho_{b'w}(\whs)\in W_k$, it follows $P\rho_{b'w}(\whs)\in W_{k-1}$, whence \[\rho_{bu}P\rho_{b'w}(\whs)\in W_s\subseteq V_0\] as $|\delta(bu)|\leq s-k+1$.  Applying the Standard Argument results in a word $t$ with $|St\inv|>|S|$ and
\begin{align*}
|t|&\leq |bu|+c_k+|b'w|\\ &\leq 1+\dim W_s-\dim W_k+c_k+ 1+\dim W_{k-1}-1 \\ &\leq n-(\dim W_k-\dim W_{k-1})+c_k= n.
\end{align*}
 This completes the proof.
\end{proof}

\subsection{Products of cyclic groups of prime order}
Let $p$ be a prime and $m\geq 1$.   Consider the group $G=\ZZ_p^m$.  Then every generating set for $G$ contains a basis and so to prove that $G$ is a {\v{C}}ern{\'y} group, it suffices to show that the Cayley graph of $G$ with respect to a basis is a {\v{C}}ern{\'y} Cayley graph.

Let's first describe the irreducible representations of $G$.  We have already seen the irreducible representation $\p\colon \ZZ_p\to \End \QQ {\QQ(\omega_p)}$ which sends the generator to left multiplication by $\omega_p$. Hence if $\psi\colon \ZZ_p^m\to \ZZ_p$ is any non-zero (and hence onto) linear functional, then the composition $\p \psi$ is an irreducible representation of $\ZZ_p^m$.  Now if $\chi$ is the character of $\p$, then $\chi\psi$ is the character of $\p\psi$.  A straightforward computation yields \[\chi(k) = \begin{cases} p-1 & k=0\\ -1 & k\neq 0\end{cases}\] (since $\chi$ summed with the trivial character of $\ZZ_p$ gives the regular representation of $\ZZ_p$).  Thus if $\psi_1,\psi_2$ are two non-zero linear functionals, then $\chi\psi_1=\chi\psi_2$ if and only if $\ker \psi_1=\ker \psi_2$.  But two non-zero functionals on a finite dimensional vector space have the same hyperplane as a kernel if and only if they are scalar multiples of each other.  In particular, the number of isomorphism classes of  irreducible representations of $\ZZ_p^m$ of the form $\p\psi$ with $\psi$ a non-zero functional equals the number of lines in the dual vector space of $\ZZ_p^m$, which is of course $(p^m-1)/(p-1)$.

Thus we have found $(p^m-1)/(p-1)$ pairwise non-isomorphic irreducible representations of degree $p-1$.  The direct sum of all these representations and the trivial representation gives a subrepresentation of the regular representation of $G$ of degree $p^m$ and so it must be the regular representation.  Thus the above representations, along with the trivial representation, constitute all the irreducible representations of $G$. Consequently, $m(G)=p-1$ while the diameter of the Cayley graph is $m(p-1)$.  In particular, for $m>1$, Theorem~\ref{mainthm} does not help us prove that $G$ is a {\v{C}}ern{\'y} group.
Nonetheless, we can show that $\ZZ_p^m$ is a {\v{C}}ern{\'y} group for all $m$.

\begin{Thm}\label{pxp}
Let $p$ be a prime.  Then $\ZZ_p^m$ is a {\v{C}}ern{\'y} group for all $m\geq 1$.
\end{Thm}
\begin{proof}
Let $G=\ZZ_p^m$ and suppose $(G,\Sigma)$ is a synchronizing automaton with $\Sigma$ containing a basis $\Delta$ for $G$.  Set $n=p^m$.
Let $S$ be a subset of $G$ with $2\leq |S|<n$.  We show that there is a word $t\in\Sigma^*$ of length at most $n$ with $|St\inv|>|S|$.  Let us assume the Standard Setup.

Since $W_0$ is $G$-invariant, we may write it as $M_1\oplus\cdots \oplus M_k$ where the subspaces $M_1,\ldots, M_k$ carry non-trivial irreducible subrepresentations of $G$.  Then there exist non-zero linear functionals $\psi_1,\ldots,\psi_k$ on $\ZZ_p^m$ so that $M_i$ affords a representation isomorphic to $\p\psi_i$ with $\p$ as in the discussion preceding the proof.  In particular, $c_0=\dim W_0=k(p-1)$.  The representation afforded by $W_0$ is $\psi=\p\psi_1\oplus\cdots \oplus \p\psi_k$ and hence, since $\p$ is injective, $\ker \psi = \bigcap_{i=1}^k \ker \psi_i$.  But $G/\ker \psi_i\cong \ZZ_p$, for all $i=1,\ldots,k$,  so $H=G/\ker\psi$  is isomorphic to a subgroup of $\ZZ_p^k$ and hence has dimension at most $k$ as a $\ZZ_p$-vector space.  Since $\Delta$ is a basis for $G$, the image of $\Delta$ is a spanning set for $H$ and hence some subset of $\Delta$ of size at most $k$ maps to a basis of $H$.  Thus each coset of $H$ can be represented by an element of $\Delta^*$ of length at most $k(p-1)= c_0$.  An application of Lemma~\ref{killwithcyclic} (with $A=G$) provides the desired word $t$.
\end{proof}

\begin{Rmk}
Notice that Theorem~\ref{pxp} only uses the case of Lemma~\ref{killwithcyclic} where $k=0$, which is the easier case.
\end{Rmk}

Using similar techniques it can also be shown that if $p_1,\ldots, p_k$ are distinct odd primes, then the Cayley graph of $G=\ZZ_{p_1}^{m_1}\times \cdots\times \ZZ_{p_k}^{m_k}$ with respect to a generating set $\Delta = \bigcup_{i=1}^k\Delta_i$, where $\Delta_i$ is a basis for $\ZZ_{p_i}^{m_i}$, is a {\v{C}}ern{\'y} Cayley graph.  Here one must use that the irreducible representations of $G$ are obtained by projecting to $\ZZ_d$ where $d\mid p_1\cdots p_k$ and then acting on $\QQ(\omega_d)$.

\subsection{Affine groups}
Fix an odd prime $p$.  Then $\ZZ_p^*$ acts naturally on $\ZZ_p$ by left multiplication and we can form the semidirect product $\ZZ_p\rtimes \ZZ_p^*$, which can be identified with the affine group $AG(1,p)$ of all maps $\ZZ_p\to \ZZ_p$ of the form $x\mapsto sx+r$ with $s\in \ZZ_p^*$ and $r\in \ZZ_p$.  Now fix a subgroup $K\leq \ZZ_p^*$ and set $G=\ZZ_p\rtimes K$.  For example, the case $K=\{\pm 1\}$ results in the dihedral group $D_p$.  Put $k=|K|$.  Suppose that $\Delta$ is a generating set for $G$ so that every translation $x\mapsto x+r$ can be represented by a word over $\Delta$ of length at most $p-1$, e.g.\ if $\Delta$ contains a non-trivial translation.  Our goal is to show that $(G,\Delta)$ is a {\v{C}}ern{\'y} Cayley graph.  First let us estimate the diameter.  Denote by $A$ the normal subgroup of translations (so $A\cong \ZZ_p$).  Since $G/A \cong K$ has size $k$, it follows that each coset of $A$ has a representative of length at most $k-1$.  Since $G=\bigcup Ag$ where $g$ runs over any given set of coset representatives, we conclude that the diameter of $(G,\Delta)$ is at most $p-1+k-1=p+k-2$ by our assumption on $\Delta$.

Define a map $\p\colon G\to \End {\QQ} {\QQ(\omega_p)}$ on the basis by $\p_{(r,s)}(\omega_p^t) = \omega_p^{st+r}$ for $0\leq t\leq p-1$.  So the factor $\ZZ_p$ acts in the way to which we are already accustomed while $K$ acts via the identification of $\ZZ_p^*$ with the Galois group $\mathrm{Gal}(\QQ(\omega_p),\QQ)$.  It is routine to verify that $\p$ is a representation.  Also if $\lambda\colon K\to \End {\QQ}{\QQ^K}$ is the regular representation of $K$ and $\pi\colon G\to K$ is the projection, then $\lambda\pi\colon G\to \End {\QQ}{\QQ^K}$ is a representation.

\begin{Prop}\label{decomp}
The regular representation of $G$ over $\QQ$ decomposes as $\lambda\pi\oplus k\cdot \p$.
\end{Prop}
\begin{proof}
We compare characters.  Let $\chi$ be the character of the regular representation of $G$. It is well known and easy to see that \[\chi(r,s) = \begin{cases}|G| & (r,s)=(0,1)\\ 0 & \text{otherwise.}\end{cases}\]  Let $\theta$ be the character of  $\lambda\pi$ and $\zeta$ the character of $\p$.  Then we have \[\theta(r,s) = \begin{cases}k & s=1\\ 0 & s\neq 1\end{cases}\]

To compute $\zeta$, first let $\alpha$ be the character of the representation $\psi$ of $G$ on $\QQ[x]/(x^p-1)$ given by $\psi_{(r,s)}(x^t+(x^p-1))= x^{st+r}+(x^p-1)$.  Then as a representation of $G$, $\QQ[x]/(x^p-1)$ decomposes as the direct sum $\QQ\oplus \QQ(\omega_p)$ where the factor $\QQ$ is spanned by $1+x+\cdots+x^{p-1}+(x^p-1)$, which is fixed by $G$ (since $G$ is a group of permutations of $\ZZ_p$ and the latter can be identified with the cyclic group $\langle x+(x^p-1)\rangle$).  Thus $\psi$ is the direct sum of $\p$ and the trivial representation.  Now $\alpha$ counts the number of $0\leq t\leq p-1$ so that $st+r\equiv t\bmod p$.  But this latter congruence is equivalent to $t(1-s)\equiv r\bmod p$ and so has $p$ solutions if $r=0,s=1$, no solutions if $s=1,r\neq 0$ and one solution otherwise.  Since $\zeta(r,s) = \alpha(r,s)-1$,  it follows \[\zeta(r,s) = \begin{cases} p-1 & r=0,s=1\\ -1 & r\neq 0,s=1\\ 0 & \text{else.}\end{cases}\]

Putting it all together, we compute \[(\theta+k\cdot \zeta)(r,s) = \begin{cases}k+k(p-1)=kp & r=0,s=1\\ k-k=0 & r\neq 0, s=1\\ 0 & \text{else} \end{cases}\] and so $\chi = \theta+k\cdot \zeta$, completing the proof.
\end{proof}

The proposition immediately leads us to deduce that $m(G)=p-1$ and consequently Theorem~\ref{mainthm} is to weak to establish that $(G,\Delta)$ is a {\v{C}}ern{\'y} Cayley graph.  Nonetheless, it is a {\v{C}}ern{\'y} Cayley graph as the following result shows.

\begin{Thm}\label{affine}
Let $K\leq \ZZ_p^*$ be a subgroup with $p$ an odd prime.  Set $G$ equal to the semidirect product $\ZZ_p^*\rtimes K$, which we view as a subgroup of the affine group $AG(1,p)$. Let $\Delta$ be a generating set for $G$ so that each translation has a representative in $\Delta^*$ of length at most $p-1$.   Then the Cayley graph $(G,\Delta)$ of $\ZZ_p\rtimes K$ is a {\v{C}}ern{\'y} Cayley graph.
\end{Thm}
\begin{proof}
If $K$ is trivial, then there is nothing to prove since we already know $\ZZ_p$ is a {\v{C}}ern{\'y} group.  So assume $K\neq 1$.
We retain the notation above and assume the Standard Setup.  We must find $t\in \Sigma^*$ with $|St\inv|>|S|$ and $|t|\leq n$.  Recalling that we have show under the hypotheses of the theorem that $\diam \Delta G\leq p-1+k-1$, if $c_r\geq 2(p-1)$ some $0\leq r\leq s$, then the Gap Bound provides a word $t$ with $|St\inv|>|S|$ and \[|t|\leq n-c_r+\diam \Delta G\leq n-2(p-1)+p-1+k-1\leq n.\] If $\dim W_s\leq n-1-2(p-1)$, then the Gap Bound again asserts the existence of a word $t$ of length no more than $n-2(p-1)+\diam \Delta G\leq n$ so that $|St\inv|>|S|$.

Next suppose that $c_r=p-1$ for some $0\leq r\leq s$.  Since $V\cong \QQ^K\oplus k\cdot \QQ(\omega_p)$ and $\dim \QQ^K/V_0=k-1<p-1$, it must be the case that $W_r=W_{r-1}\oplus U_r$ with $U_r\cong \QQ(\omega_p)$ (where we take $W_{-1}=0$, as usual).  But if $A$ is the subgroup of translations, then $U_r$ affords a non-trivial irreducible representation of $A$, whence Lemma~\ref{killwithcyclic} provides the desired word $t$ as by assumption each element of $A$ has a representative in $\Delta^*$ of length at most $p-1$ and $c_r=p-1$.

If we are in none of the above cases, then $W_s$ must contain as constituents at least $k-1$ of the $k$ copies of $\QQ(\omega_p)$. In the notation of the Standard Setup,  $W_s$ decomposes as  $U_0\oplus U_1\oplus\cdots\oplus U_s$ with $\dim U_r = c_r$.  Here no $U_i\cong \QQ(\omega_p)$ or contains $\QQ(\omega_p)$ as a constituent with multiplicity greater than $1$, or we would be back in one of the previous cases.
From the fact that $\QQ^K/V_0$ has at most $k-1$ irreducible constituents, it follows that $s=k-2$ and each $U_r\cong \QQ(\omega_p)\oplus M_r$ where $M_r$ is a non-trivial irreducible constituent of $\QQ^K$, for $0\leq r\leq s$.  Thus $V_0\cong W_{s}\oplus \QQ(\omega_p)$ and hence $\dim W_s\leq n-1-(p-1)$ from which there results, by the Gap Bound, a word $t$ with $|St\inv|>|S|$ and length at most \[1+\dim W_s-(p-1)+\diam \Delta G\leq n-2(p-1)+p-1+k-1\leq n.\] This completes the proof, establishing the theorem.
\end{proof}

\begin{Rmk}
Let us remark that the last case of the above proof can only happen when $k=2$ since if $K$ has $k-1$ non-trivial irreducible representations, then each of them must have degree $1$ and so $m(K)=1$, which implies $K\cong \ZZ_2^m$.  But $K$ must be cyclic, being a subgroup of $\ZZ_p^*$, and consequently $k=2$, as claimed.
\end{Rmk}

An important special case of Theorem~\ref{affine} is the full affine group.

\begin{Cor}
If $p$ is an odd prime, any Cayley graph of the affine group $AG(1,p)$ with respect to a generating set containing a translation is a {\v{C}}ern{\'y} Cayley graph.
\end{Cor}

\subsection{Dihedral groups: revisited}
In this section we show that if $p$ is an odd prime, then the dihedral groups $D_p$ and $D_{p^2}$ are {\v{C}}ern{\'y} groups.  Let us begin with $D_p$.  Since the subgroup of rotations of a regular $p$-gon is cyclic of prime order, and hence generated by any non-trivial element, there are two types of generating sets for $D_p$ that are minimal with respect to containment:  either a reflection and a rotation, or two distinct reflections.  Indeed, any generating set $\Delta$ must contain a reflection.  If $\Delta$ contains a rotation, then we are in the first case; if $s_1,s_2\in \Delta$ are distinct reflections, then $s_1s_2$ is a rotation by twice the angle between their respective lines of reflection and hence $s_1,s_2$ generates the dihedral group.

A similar analysis holds for $D_{p^2}$.  Let $r$ be a rotation of order $p^2$ and let $K$ be the subgroup generated by $r^p$.  Then $K$ is a normal subgroup and $D_{p^2}/K\cong D_p$.   We claim that $\Delta$ is a generating set of for $D_{p^2}$ if and only if under the canonical projection $\rho\colon D_{p^2}\to D_p$ one has that $\rho(\Delta)$ generates $D_p$.  Necessity is clear.  For sufficiency, observe that if $\rho(\Delta)$ is a generating set, then either it contains a reflection and a rotation or two reflections.  Consider the first case.  Then the rotation is of the form $aK$ where $a$ is a rotation not belonging to $K$.  But any element of $\langle r\rangle$ not belonging to $K$ is a generator.  Thus $a$ is a rotation of order $p^2$ and $\Delta$ generates $D_{p^2}$.  In the second case, we have reflections $s_1,s_2$ so that $s_1K,s_2K$ generate $D_p$.  Then $s_1s_2K$ is a non-trivial rotation and so $s_1s_2\notin K$.  Hence,  $s_1s_2$ generates $\langle r\rangle$ and so $s_1,s_2$ generate $D_{p^2}$.  It follows that minimal generating sets of $D_{p^2}$ with respect to containment consist either of a reflection and rotation of order $p^2$ or of two reflections $s_1,s_2$ so that $s_1s_2$ is a rotation of order $p^2$.  The reader should note that the same argument applies \textit{mutatis mutandis} to $D_{p^m}$.

In Subsection~\ref{ss:dihedralgroups}, we showed that Cayley graphs of $D_p$ and $D_{p^2}$ with respect to a generating set consisting of a rotation and a reflection are {\v{C}}ern{\'y} Cayley graphs (the former is also covered by Theorem~\ref{affine}), so we are left with considering generating sets consisting of two reflections.

Consider for the moment $D_n$ with $n$ odd.  Let $s,s'$ be two reflections so that $ss'$ is a rotation of order $n$.  Then we claim that the diameter of $D_n$ is at most $n$ (actually it is exactly $n$, as is well known in the theory of reflection groups).  Indeed, since $s,s'$ are involutions, it follows that $(ss')\inv =s's$ and so each non-trivial rotation can be written uniquely in the form $(ss')^k$ or $(s's)^k$ with $0\leq k\leq \frac{n-1}{2}$, that is each rotation can be represented by a word of length at most $n-1$.  Since each reflection is a product of $s$ with a rotation, this gives the upper bound of $n$.  It now follows that Theorem~\ref{affine} applies to show that $D_p$ is a {\v{C}}ern{\'y} group.

\begin{Thm}
Let $p$ be an odd prime.  Then the dihedral group $D_p$ of order $2p$ is a {\v{C}}ern{\'y} group.
\end{Thm}
\begin{proof}
Viewing $D_p$ as a subgroup of the affine group $AG(1,p) = \ZZ_p\rtimes \ZZ_p^*$, the above discussion shows that each translation can be represented by a word of length at most $p-1$ for any generating set of $D_p$.  Theorem~\ref{affine} then provides the desired conclusion.
\end{proof}

To prove that $D_{p^2}$ with $p$ an odd prime is a {\v{C}}ern{\'y} group we first need to decompose the regular representation of $D_{p^2}$ over the rational numbers.   Let $r$ be a rotation by $2\pi/p^2$ and $s$ a reflection over an axis of symmetry of the regular $p^2$-gon.  Let $\alpha\colon D_{p^2}\to \QQ^*$ be given by sending each reflection to $-1$ and rotation to $1$.  Also note that $\QQ(\omega_p)$ and $\QQ(\omega_{p^2})$ afford irreducible representations of $D_{p^2}$ by having $r$ act as multiplication by $\omega_p,\omega_{p^2}$, respectively, and $s$ acting as complex conjugation.  Again the latter two representations are already irreducible when restricted to $\langle r\rangle$.

\begin{Prop}\label{regrepofDp2}
Let $p$ be an odd prime.  Then the regular representation of $D_{p^2}$ decomposes as the direct sum of the trivial representation, $\alpha$ and two copies of both $\QQ(\omega_p)$ and $\QQ(\omega_{p^2})$.
\end{Prop}
\begin{proof}
For notational purposes let $r$ be a rotation by $2\pi/p^2$ and $s$ a reflection.
Let $\chi_1,\chi_2$ be the characters afforded by $\QQ(\omega_p)$ and $\QQ(\omega_{p^2})$ respectively.  Notice that $\alpha$ can be viewed as its own character.  We show that the character $\chi$ of the regular representation is the sum of the trivial character $\tau$ with $\alpha+2\cdot \chi_1+2\cdot \chi_2$.
Since the value of a character at $1$ is its degree, first note
\begin{align*}
\tau(1)+\alpha(1)+2\chi_1(1)+2\chi_2(1)&= 1+1+2\phi(p)+2\phi(p^2)\\ &= 1+1+2(p-1)+2(p^2-p) = 2p^2\\ &=\chi(1).
\end{align*}
Next we remark that $\chi(g)=0$ all $1\neq g\in D_{p^2}$.    From the computation in Proposition~\ref{decomp} for $\zeta$, it follows that \[\chi_1(r^k) = \begin{cases} p-1 & p\mid k\\ -1 & p\nmid k\end{cases}\] while $\chi_1(sr^k) = 0$ all $k$.

Since the regular representation of $\ZZ_{p^2}$ is $\QQ\oplus \QQ(\omega_p)\oplus \QQ(\omega_{p^2})$ we may deduce that \[\chi_2(r^k) = \begin{cases} -1-(p-1) = -p & p\mid k, k\neq 0\\ -1-(-1) = 0 & p\nmid k.\end{cases}\] We claim that $\chi_2(sr^k)=0$ all $k$.  Since all the reflections are conjugate in $D_n$ with $n$ odd (rotation acts transitively on the axes of symmetry of a regular $n$-gon with $n$ odd), it suffices to deal with $\chi_2(s)$ (recall characters are traces and similar linear operators have the same trace).

To ease notation, set $\omega =\omega_{p^2}$.  Then $\{1,\omega,\ldots,\omega^{p^2-p-1}\}$ is a basis for $\QQ(\omega_{p^2})$ and the minimal polynomial for $\omega$ is the cyclotomic polynomial \[1+x^p+(x^p)^2+\cdots +(x^p)^{p-1}.\]  If $p<m<p^2-p$, then $\ov{\omega^m} = \omega^{p^2-m}$ and $p^2-m<p^2-p$.  Since $p^2-m\neq m$, we conclude basis vectors of this form do not contribute to the trace of the operator complex conjugation.  From the minimal polynomial for $\omega$ it follows \[\ov{\omega^p} = \omega^{p^2-p} = -1-\omega^p-(\omega^p)^2-\cdots -(\omega^p)^{p-2}\] and so the basis vector $\omega^p$ contributes $-1$ to the trace of complex conjugation as an operator.  If $0<m<p$, then \[\ov{\omega^m} = \omega^{p^2-m} = \omega^{p^2-p}\omega^{p-m} = (-1-\omega^p-(\omega^p)^2-\cdots -(\omega^p)^{p-2})\omega^{p-m}.\] Note that $kp+p-m = m$ with $0\leq k\leq p-2$ implies $(k+1)p=2m$, a contradiction since $2,m<p$.  So basis vectors of this form do not contribute to the trace.  Finally, $\ov 1=1$ and so the basis vector $1$ contributes $1$ to the trace.  Thus the trace of complex conjugation is zero, i.e.\ $\chi_2(s)=0$, as was required.

It follows $\tau(sr^k)+\alpha(sr^k)+2\chi_1(sr^k)+2\chi_2(sr^k) = 1-1+0+0 = 0$ and
\begin{align*}
\tau(r^k)+\alpha(r^k)+2\chi_1(r^k)+2\chi_2(r^k)&=\begin{cases} 1+1+2(p-1)-2p  & p\mid k, k\neq 0\\
1+1+2(-1)+0 &p\nmid k\end{cases} \\ &=0
\end{align*}
establishing the desired equality $\chi = \tau+\alpha+2\cdot \chi_1+2\cdot \chi_2$.
\end{proof}

\begin{Thm}
Let $p$ be an odd prime.  Then the dihedral group $D_{p^2}$ of order $2p^2$ is a {\v{C}}ern{\'y} group.
\end{Thm}
\begin{proof}
By the discussion at the beginning of this subsection we need only handle the case that the generating set $\Delta$ consists of two reflections $s,s'$ with $r=ss'$ a reflection of order $p^2$.  Let us assume the Standard Setup and prove the existence of a word $t$ of length at most  $n=2p^2$ so that $|St\inv|>|S|$.

As shown above, $\diam \Delta {D_{p^2}}\leq p^2$.  Also $V_0$ has five irreducible constituents: $\alpha$ of degree $1$, two copies of $\QQ(\omega_p)$ each of degree $p-1$ and two copies of $\QQ(\omega_{p^2})$ each of degree $p^2-p$.   Let $K=\langle r^p\rangle$; so $D_{p^2}/K \cong D_p$ and $sK,s'K$ generate the quotient group. Notice that $K$ is the kernel of the representation of $D_{p^2}$ on $\QQ(\omega_p)$.

Recalling $W_s=U_0\oplus\cdots\oplus U_s$,
assume first that $U_i$ affords $\alpha$ for some $0\leq i\leq s$.  Then applying Lemma~\ref{killwithcyclic} to $A=\langle s\rangle$ establishes the existence of the desired word $t$.

Next assume that $U_i\cong \QQ(\omega_p)$ for some $0\leq i\leq s$.   Let $A=\langle r\rangle$ be the subgroup of rotations.  Then $\QQ(\omega_p)$ affords a non-trivial irreducible representation $\psi$ of $A$ and $\ker \psi = K$.  Since every rotation in $D_p\cong D_{p^2}/K$ can be written as either $(ss'K)^m$ or $(s'sK)^m$ with $m\leq \frac{p-1}{2}$, it follows that each coset of $A/K$ has a representative from $\Delta$ of length at most $p-1$ and so Lemma~\ref{killwithcyclic} again applies to guarantee the desired word $t$ exists.

Suppose that, for some $0\leq i\leq s$, we have $U_i$ is isomorphic to either $2\QQ(\omega_p)$, $\alpha\oplus\QQ(\omega_p)$ or $\alpha\oplus 2\QQ(\omega_p)$.  Let $\psi\colon G\to \End {\QQ} {U_i}$ be the representation afforded by $U_i$.  Then $\ker U_i = K$ and $c_i\geq p\geq \diam {sK,s'K} {D_{p^2}/K}$ and so an application of Lemma~\ref{killwithcyclic} yields the sought after word $t$.

We claim that in all other cases, the Gap Bound provides the desired conclusion.  First we claim that unless there exists $0\leq i<j\leq s$ so that $U_i$ and $U_j$ both have $\QQ(\omega_{p^2})$ as constituents, the Gap Bound immediately provides the result.  Indeed, if no copy of $\QQ(\omega_{p^2})$ is a constituent of $W_s$, then \[1+\dim W_s+\diam \Delta {D_{p^2}}\leq n- 2(p^2-p)+p^2\leq n\] and the Gap Bound establishes the desired result.   On the other hand, if exactly one copy of  $\QQ(\omega_{p^2})$ is a constituent of $W_s$, then $c_r\geq p^2-p$ some $0\leq r\leq s$ and also $1+\dim W_s\leq n-(p^2-p)$.  So the Gap Bound yields a word $t$ of length at most $n-(p^2-p)-(p^2-p)+p^2\leq n$ in this case as well.

So let $U_i,U_j$ be as above.  If no constituent of $W_s$ is isomorphic to $\QQ(\omega_p)$, then again the Gap Bound provides the desired result since \[1+\dim W_s-(p^2-p)+p^2\leq n-2(p-1)-(p^2-p)+p^2 = n-p+2\leq n\] as $p\geq 3$.  If we are not in one of the cases previously considered, then $\QQ(\omega_p)$ may only occur as a constituent of $U_i$ or $U_j$ in $W_s$. If $\QQ(\omega_p)$ is a constituent of either $U_i$ or $U_j$, but not both, then the Gap Bound once again yields the desired result since  \[1+\dim W_s-(p^2-p+p-1)+p^2\leq n-(p-1)-(p^2-1)+p^2 = n-p+2\leq n.\]  Thus we are left with the case that $U_i$ and $U_j$ each have $\QQ(\omega_p)$ and $\QQ(\omega_{p^2})$ as constituents.  In particular, we have $c_i,c_j\geq p^2-1$.

Now again, by the cases previously considered, either $\alpha$ is not a constituent of $W_s$ or $\alpha$ is a constituent of $U_i$ or $U_j$.  But then again the Gap Bound handles the result since in the latter case either $c_i$ or $c_j$ is $p^2=\diam \Delta D_{p^2}$, while in the former $1+\dim W_s = n-1$ and so the Gap Bound yields $n-1-(p^2-1)+p^2=n$ as an upper bound on the length of $t$.  This completes the proof.
\end{proof}

\section{Open questions}
There are a number of open questions left by this paper.  As it is not quite clear that the {\v{C}}ern{\'y} conjecture is true --- there is not even a quadratic bound at the full level of generality --- the fact that there are quadratic bounds in the context of this paper makes the following question enticing.

\begin{Question}
Is it true that all groups are {\v{C}}ern{\'y} groups?
\end{Question}

Dubuc's work~\cite{dubuc} begs the question as to whether all cyclic groups are {\v{C}}ern{\'y} groups.

\begin{Conjecture}
All cyclic groups are {\v{C}}ern{\'y} groups.
\end{Conjecture}

The difficulty in working on this conjecture is that Dubuc seems to use in an essential way that each element of a cyclic group of order $n$ has a \emph{unique} representation by a word of length at most $n-1$ with respect to a cyclic generating set. I suspect that a little bit of number theory may be needed in the general case.

The next natural step would be to consider abelian groups.    I would guess
that if one can handle the above conjecture, then the next conjecture should be accessible.

\begin{Conjecture}
All abelian groups are {\v{C}}ern{\'y} groups.
\end{Conjecture}

I suspect that Dubuc's techniques~\cite{dubuc} can be extended to show that the Cayley graph of a dihedral group with respect to a generating set consisting of a reflection and a rotation is a {\v{C}}ern{\'y} Cayley graph.  I will put forth the following bolder conjecture.

\begin{Conjecture}
Dihedral groups are {\v{C}}ern{\'y} groups.
\end{Conjecture}

Finally, given the large degrees of representations and the substantial amount of knowledge in the literature concerning representations of symmetric groups, it seems natural to ask:

\begin{Question}
Are all symmetric groups {\v{C}}ern{\'y} groups?
\end{Question}

\bibliographystyle{abbrv}
\bibliography{standard2}

\end{document}